\documentclass[12pt]{amsart}

\usepackage{amssymb}
\usepackage{latexsym}
\usepackage{amsmath}

\def\c{{\mathbf C}}
\def\w{{\mathcal W}}
\def\l{{\mathcal L}}
\def\p{{\mathcal P}}
\def\a{{\mathcal A}}
\def\r{{\mathcal R}}
\def\s{{\mathcal S}}

\def\uu{{\mathbf U}}
\def\ss{{\mathbf S}}
\def\oo{{\mathbf O}}
\def\nn{{\mathbf N}}
\def\ww{{\mathbf W}}
\def\gg{{\mathbf G}}
\def\ll{{\mathbf L}}
\def\rr{{\mathbf R}}
\def\zz{{\mathbf Z}}
\def\qq{{\mathbf Q}}

\def\cat{{\mathcal C}}
\def\ca{{{\mathcal C}_X}}
\def\lie{\mathfrak{g}}
\def\sl{\mathfrak{sl}}

\begin{document}

\newtheorem{defi}{Definition}[section]
\newtheorem{prop}{Proposition}[section]
\newtheorem{theo}{Theorem}[section]
\newtheorem{lemm}{Lemma}[section]
\newtheorem{coro}{Corollary}[section]
\newtheorem{conj}{Conjecture}[section]
\newtheorem{prob}{Problem}[section]

\hyphenation{Wo-ro-no-wi-cz}
\hyphenation{se-mi-ring}
\hyphenation{quan-ti-za-tion}
\hyphenation{Gro-then-dieck}
\hyphenation{mo-no-i-dal}

\title
{Fusion rules for representations of compact quantum groups}
\author{Teodor Banica}
\address{
Institut de Math\'ematiques de Jussieu, Case 191, Universit\'e Paris 6, 
4 Place Jussieu, 75005 Paris}
\email{banica@math.jussieu.fr}

\maketitle

\section*{Introduction}

The compact quantum groups are objects which generalise at the
same time the compact groups, the duals of discrete groups and the
$q-$deformations (with $q>0$) of classical compact Lie groups. A
compact quantum group is an abstract object which may be
described by (is by definition the dual of) the algebra of ``continuous
functions on it'', which is a Hopf ${\bf C}^*$-algebra. A system of
axioms for Hopf ${\bf C}^*$-algebras which leads to a satisfactory
theory of compact quantum groups (e.g. a theorem stating the existence
of the Haar measure) was found by Woronowicz at the end of the 80's.

The representation theory of compact quantum groups gives rise to rich
combinatorial structures. By Woronowicz's analogue of the Peter-Weyl
theory, each (finite dimensional unitary) representation of a compact
quantum group $G$ is completely reducible. In particular given two
irreducible representations $a$ and $b$, their tensor product
decomposes in a unique way (up to equivalence) as a sum of irreducible representations
$$a\otimes b\simeq c+d+e+\cdots$$

These formulae are called fusion rules for irreducible representations
of $G$. The fusion semiring $\r^+(G)$ is by definition the set of
equivalence classes of finite dimensional continuous representations
of $G$, endowed with the binary operations $+$ (the sum of classes of
corepresentations) and $\otimes$ (the tensor product of classes of
corepresentations). It is the algebraic structure
describing the collection of all fusion rules. There are two basic examples: 

\noindent -- if $G$ is a compact group then $\r^+(G)$ is the usual fusion semiring of $G$.

\noindent -- if $G$ is the dual of a discrete group $\Gamma$ then $\r^+(G)$ is the convolution semiring of $\Gamma$.

More generally, we have the following construction of fusion rules and semirings: in a semisimple monoidal category $\cat$,
formulae of the form $a\otimes b\simeq c+d+e+\cdots$ with
$a,b,c,d,e,\ldots$ simple objects of $\cat$, may be called fusion rules for simple objects of
$\cat$. The algebraic structure describing the collection of all
fusion rules is the Grothendieck semiring $([\cat ],+,\otimes )$ of
$\cat$. Fusion rules -- and related algebraic objects, such as fusion
semirings, rings, algebras and principal graphs -- appear in this way in many
recent theories arising from mathematics and physics, such as
conformal field theory, subfactors, quantum groups at roots of
unity. They can be thought of as being a common language for these
theories. The above fusion semiring $\r^+(G)$ arises also in this way:
it is isomorphic to the Grothendieck semiring of the semisimple
monoidal category $Rep(G)$ of finite dimensional continuous
representations of $G$.

In this paper we give a survey of some recent results on the fusion
semirings of compact quantum groups (computations of and applications
to discrete quantum groups) by using the following simplifying
terminology: we say that a compact quantum group $G$ is an $\r^+$-deformation of a compact
quantum group $H$ if their fusion semirings are isomorphic. The paper
contains also some easy related results (with proofs), two conjectures
and many remarks and comments, some of them concerning classification
by invariants related to $\r^+$.

I would like to thank Shuzhou Wang for some useful comments on a
preliminary version of this paper.

\tableofcontents

\section{The fusion semiring of a Woronowicz algebra}

The Woronowicz algebras are the Hopf $\c^*$-algebras which correspond
to both notions of ``algebras of continuous functions on compact
quantum groups'' and ``$\c^*$-algebras of discrete quantum
groups''. That is, specialists agree that one can define the category of
compact quantum groups to be the dual $\widehat{\w}$ of the category
$\w$ of Woronowicz algebras, and the category of discrete quantum groups to be $\widehat{\widehat{\w}}=\w$.

The category $\w$ may be defined in the following two equivalent ways.

-- the ``algebras of continuous functions on compact matrix quantum groups''
or ``$\c^*$-algebras of finitely generated discrete quantum
groups'' were introduced by Woronowicz in \cite{wo2} via a fairly
simple and useful set of axioms (see definition 1.1 below; see also \cite{wo3}, \cite{wo4} and section 3
below for various non-trivial translations of these axioms). Such objects are
called finitely generated Woronowicz algebras. With
a suitable definition for morphisms one gets a category $\w_{f.g.}$,
and $\w$ may be defined as $\w =Ind(\w_{f.g.})$, the category of
inductive limits of $\w_{f.g.}$. See e.g. Baaj-Skandalis \cite{bs}.

-- $\w$ is the category of bisimplifiable unital Hopf
$\c^*$-algebras. See Woronowicz \cite{wo5}.

Almost all known examples of Woronowicz algebras (see sections 2
and 8 below for a quite complete list) are finitely generated. Also, the
first approach to $\w$ seems to be more useful than the second one. In
what follows we will mainly be interested in objects of $\w_{f.g.}$
and we will always refer the reader to the fundamental paper
\cite{wo2}. We mention that all basic results are known to extend easily from $\w_{f.g.}$ to $\w$ (see e.g. \cite{bs}, \cite{wo5}).

\begin{defi}[cf. \cite{wo2}]
A finitely generated (or co-matricial) Woronowicz algebra is a pair
$(A,u)$ consisting of a unital $\c^*$-algebra $A$ and a unitary matrix
$u\in M_n(A)$ subject to the following conditions:

(i) the coefficients of $u$ generate in $A$ a dense $*$-subalgebra, denoted $\a$.

(ii) there exists a $\c^*$-morphism $\delta :A\rightarrow A\otimes_{min} A$ such that $(id\otimes\delta )u=u_{12}u_{13}$.

(iii) there exists a linear antimultiplicative map $\kappa
:\a\rightarrow\a$ such that $\kappa *\kappa *=id$ and such that 
$(id\otimes\kappa )u=u^{-1}$. 
\end{defi}

The dense subalgebra $\a$ is an (involutive) Hopf $\c$-algebra with
the restriction of $\delta$ as comultiplication, with counit defined
by $\varepsilon :u_{ij}\mapsto\delta_{i,j}$, and with antipode
$\kappa$ (see \cite{wo2}). We recall that if $V$ be a finite
dimensional $\c$-linear space, a coaction of the Hopf $\c$-algebra $\a$ on $V$ is a linear map $\beta :V\to V\otimes\a$ satisfying
$$(id\otimes\delta )\beta =(\beta\otimes id )\beta ,\,\,\,\,
(id\otimes\varepsilon )\beta =id$$
A corepresentation of $\a$ on $V$ is an element $u\in\l (V)\otimes\a$
satisfying
$$(id\otimes\delta )u=u_{12}u_{13},\,\,\,\, (id\otimes\varepsilon )u=1$$
Coactions and corepresentations are in an obvious one-to-one correspondence. We
prefer to work with corepresentations. We recall that if $v,w$
are corepresentations of $\a$ on $V,W$ then their sum $v+w$ is the
corepresentation $diag(v,w)$ of $\a$ on $V\oplus W$, and their tensor
product $v\otimes w$ is the
corepresentation $v_{13}w_{23}$ of $\a$ on $V\otimes W$. $v$ and $w$ are said to be equivalent if there exists an
invertible linear map $T\in \l (V,W)$ such that $w=(T\otimes
1)v(T^{-1}\otimes 1)$.

The finite dimensional corepresentations of the Hopf $\c$-algebra $\a$
will be called finite dimensional smooth corepresentations of the finitely
generated Woronowicz algebra $(A,u)$. We will sometimes call them ``finite
dimensional corepresentations'', or just ``corepresentations'' (there
will be no other kind of corepresentation to appear in this paper).

Given a finitely generated Woronowicz algebra $(A,u)$ one can
construct its ``full version'' $(A_{full},u)$ and its ``reduced
version'' $(A_{red},u)$, which are in general different from $(A,u)$. This is because the Haar functional is not necessarely
faithful (see \cite{wo2}, \cite{bs}). There exist canonical morphisms
of $\c^*$-algebras
$$A_{full}\to A\to A_{red}$$
$A$ is said to be full if $A_{full}\to A$ is an isomorphism, reduced if
$A\to A_{red}$ is an isomorphism and amenable (as a Woronowicz algebra) if $A_{full}\to A_{red}$ is an isomorphism.

The morphisms between two finitely generated Woronowicz algebras
$(A,u)$ and $(B,v)$ are by definition the morphisms of Hopf
$*$-algebras $\a\to {\mathcal B}$. Each such morphism canonically
extends to a $\c^*$-morphism $A_{full}\to B_{full}$. Notice that:

-- $(A,u)$ is isomorphic to both $(A_{full},u)$ and $(A_{red},u)$.

-- if $v$ is a finite dimensional corepresentation of $\a$ whose
coefficients generate $\a$ as a $*$-algebra then $(A,u)$ is isomorphic
to $(A,v)$.

This definition of morphisms is the one which makes the category of
finitely generated discrete groups to embed in a unique (covariant) way
into the category of finitely generated Woronowicz algebras. Indeed,
if $\Gamma$ is a finitely generated discrete group then
$$(A,u)=(\c^*_\pi (\Gamma ), diag(\pi (g_i)))$$
is a finitely generated Woronowicz
algebra, for any set $\{ g_i\}$ of generators of $\Gamma$ and for any
unitary faithful representation $\pi$ of $\Gamma$ such that
$\pi\otimes\pi$ is contained in a multiple of $\pi$ (see \cite{wo2}).

To simplify the terminology we will sometimes call ``finitely
generated Woronowicz algebra'' the isomorphism class of a finitely
generated Woronowicz algebra, and ``corepresentation'' the equivalence
class of a corepresentation. The places where this abuse of language will be used should be clear to the reader.

These notions extend to all Woronowicz algebras (see e.g. \cite{bs}, \cite{wo5}).

\begin{defi}
The fusion semiring $R^+(A)$ of a Woronowicz algebra $A$ is the set of
equivalence classes of finite dimensional smooth corepresentations of
$A$, endowed with the binary operations $+$ (the sum of classes of
corepresentations) and $\otimes$ (the tensor product of classes of corepresentations).
\end{defi}

By cosemisimplicity coming from Woronowicz's Peter-Weyl type theory
\cite{wo2} $R^+(A)$ is isomorphic as an additive monoid to the free monoid $\nn\cdot Irr(A)$, where $Irr(A)$ is the set of
equivalence classes of finite dimensional irreducible corepresentations of
$A$.
$$(R^+(A),+)\simeq \nn\cdot Irr(A)$$
Thus $R^+(A)$ encodes the same information as the collection of formulae of the form
$$a\otimes b=c+d+e+\cdots$$
with $a,b,c,d,e,\ldots\in Irr(A)$. These formulae, which describe the
splitting of a tensor product of irreducible corepresentations into a sum of irreducible corepresentations, are called fusion rules for
irreducible corepresentations of $A$.

Notice that the complex conjugation of corepresentations
$u\mapsto\overline{u}$ makes $R^+(A)$ an involutive semiring, but the involution {\em is not} an extra structure of the semiring $R^+(A)$: if
$u\in Irr(A)$ then $\overline{u}$ may be characterized as being the
unique $v\in Irr(A)$ such that $1\subset u\otimes v$, and then -- may be extended by additivity to the whole $R^+(A)$.

A few remarks on some related algebraic objects, that will not appear
in the rest of the paper. The fusion (or representation)
ring is the Grothedieck ring
$$R(A)=K(R^+(A))$$
By extending
scalars we get the fusion $\c$-algebra
$$\c Alg(A)=R(A)\otimes_\zz\c$$
The involution of $R^+(A)$ extends by antilinearity to an involution
of $\c Alg(A)$. The linear map $\phi\in\c Alg(A)^*$ given by $\phi (a)=\delta_{a,1}$ for all $a\in
Irr(A)$ is a faithful positive trace. By applying the GNS
construction to $(\c Alg(A),\phi )$ we get the fusion $\c^*$-algebra
$$\c^* Alg(A)=\overline{\pi_\phi (\c Alg(A))}^{\mid\mid\ \mid\mid}$$
By \cite{wo2} $\c Alg(A)$ embeds into $\a$ via the linear
extension of the character application
$$u\mapsto\chi (u)=(Tr\otimes id)u$$
The linear
form $\phi$ corresponds in this way to the restriction of the Haar
functional $\int :\a\to\c$, so the embedding $\chi$ extends to an
embedding $\c^* Alg(A)\to A_{red}$, also denoted by $\chi$. The image of $\c^* Alg(A)$ is the
$\c^*$-algebra of ``central functions on the corresponding compact
quantum group''
$$A_{central}=\chi (\c^* Alg(A))$$
Each isomorphism of the form $R^+(A)\simeq
R^+(B)$ canonically extends to isomorphisms between the above four algebraic objects associated to $A$ and $B$. However, the converses are not true, i.e. the
semiring $R^+(A)$ may contain (much) more information on $A$ than
$R(A)$, $\c Alg(A)$, $\c^* Alg(A)$ or $A_{central}$. This is the reason why we
have to use the quite unfamiliar semirings.

Note that $a\geq b$ $\Longleftrightarrow$ $a-b\in R^+(A)$ makes the
fusion ring $R(A)$ an ordered ring, which contains the same
information on $A$ as $R^+(A)$. We prefer to use the semiring $R^+(A)$
instead of the ordered ring $(R(A),\geq )$.

Let us also mention that the dimension and quantum dimension of
corepresentations are morphisms of semirings
$$dim:R^+(A)\to (\nn ,+,\cdot ),\,\,\, qdim:R^+(A)\to (\rr_+^*
,+,\cdot )$$
Both of them may happen to be extra structures of $R^+(A)$
(cf. \cite{cras}, resp. \cite{wo1}). They will be used only at the end
of the paper.

In the commutative and cocommutative cases, the fusion semiring can be
computed as follows.

\begin{theo}[\cite{wo2}]
Let $G$ be a compact group. Then $A=C(G)$ is a Woronowicz
algebra. The finite dimensional corepresentations of $A$ are in one-to-one
correspondence with the finite dimensional representations of $G$, and
this induces an isomorphism between $R^+(A)$ and the usual fusion
semiring $\r^+(G)$.
\end{theo}

\begin{theo}[\cite{wo2}]
Let $\Gamma$ be a discrete group. Then $A=\c^*(\Gamma )$ is
a Woronowicz algebra. The finite dimensional irreducible
corepresentations of $A$ are all 1-dimensional, and are in one-to-one
correspondence with the elements of $\Gamma$. Their fusion corresponds
in this way to the product of $\Gamma$. That is, the fusion semiring
of $\c^*(\Gamma )$ is isomorphic to the convolution semiring $\nn\cdot\Gamma$ of $\Gamma$.
\end{theo}

\section{Examples of $R^+$-deformations}

In this section we give a survey of some recent computations of fusion
semirings of Woronowicz algebras. It turns out that (the main parts of)
these results have very short statements when using the following terminology.

\begin{defi}
Let $A$ and $B$ be two Woronowicz algebras. We say that $A$ is an
$R^+$-deformation of $B$ if there exists an isomorphism of semirings $R^+(A)\simeq R^+(B)$.
\end{defi}

Later on we will use the following related definition: we say that $A$
is a dimension-preserving $R^+$-deformation of $B$ if there exists an
isomorphism of semirings $R^+(A)\simeq R^+(B)$ which preserves the
dimensions of corepresentations.

Maybe the first requirement for a good notion of deformation is that $\c^*$-algebras of discrete groups should be rigid (:=
undeformable).

\begin{prop}
If $\Gamma$ is a discrete group then $\c^*(\Gamma )$ is
$R^+$-rigid.
\end{prop}

\begin{proof}
Let $A$ be an $R^+$-deformation of
$\c^*(\Gamma )$, and choose an isomorphism $f:\nn\cdot\Gamma\simeq
R^+(A)$ (cf. theorem 1.2). For any $g\in\Gamma$ we have $gg^{-1}=1$, and by applying $f$
we get $f(g)\otimes f(g^{-1})=1$. This shows that the corepresentation $f(g)$ is $1$-dimensional, for any $g$. Since $A$ is generated by the
coefficients of its irreducible finite dimensional corepresentations,
and since $Irr(A)=f(\Gamma )$, we get that $A$ is cocommutative, hence
isomorphic to some $\c^*(\Gamma^\prime )$
(cf. \cite{wo2}). But from $\nn\cdot\Gamma\simeq \nn\cdot\Gamma^\prime$
we get $\Gamma\simeq\Gamma^\prime$, so $A\simeq\c^*(\Gamma )$.
\end{proof}

We discuss now the relationship between $R^+$-deformations and
$q$-deformations. Let $\lie$ be a complex Lie algebra of type A,B,C,D and let $G$ be the corresponding compact connected simply-connected
Lie group. Let $q\in\c^*$ be a number which is not a root of unity, and let
$U_q\lie$ be the Drinfeld-Jimbo quantization of the universal
enveloping algebra $U\lie$ \cite{d}, \cite{ji}. It follows from Rosso \cite{r2} that if
$q>0$, the restricted dual $(U_q\lie )^\circ$ has a canonical
involution, has a $\c^*$-norm and its completion is a Woronowicz
algebra, called $C(G)_q$. The results of Lusztig and Rosso on the
$q$-deformation of finite dimensional representations of $U\lie$
\cite{l}, \cite{r1} show via \cite{r2} that these $q$-deformations are
$R^+$-deformations. More generally, it follows from work of Levendorskii and
Soibelman \cite{le}, \cite{ls}, \cite{so} that given any simple
compact Lie group $G$ and any $q>0$ one can define a $q$-deformation
$C(G)_q$ of $C(G)$, which is an $R^+$-deformation of $C(G)$.

\begin{theo}[\cite{d,ji,le,ls,l,r1,r2,so}]
$C(G)_q$ is an $R^+$-deformation of $C(G)$.
\end{theo}

We mention that for $G=\ss\uu (2)$ (resp. $G=\ss\uu (N)$ for any $N$)
the construction of $C(G)_q$ and the computation of its fusion
semiring were also done by Woronowicz in \cite{wo1}
(resp. \cite{wo3}) by using different methods.

We recall that for $n\geq 2$ and $F\in\gg\ll (n,\c )$ the $\c^*$-algebra $A_u(F)$
is defined by generators $\{ u_{ij}\}_{i,j=1,\ldots ,n}$ and the
relations making the matrices $u$ and $F\bar{u}F^{-1}$ unitaries. It
is a Woronowicz algebra. Its universal property shows that it
corresponds to both notions of ``algebra of continuous functions on
the quantum (or free) unitary group'' and ``$\c^*$-algebra of the free discrete
quantum group''. See Van Daele-Wang \cite{vdw}. By \cite{cmp} the
irreducible corepresentations of $A_u(F)$ can be labeled by the
elements of the free monoid $\nn*\nn$, $Irr(A_u(F))=\{ r_x\mid x\in\nn
*\nn\}$ such that the fusion rules are
$$r_x\otimes r_y=\sum_{x=ag,\, y=\overline{g}b}r_{ab}$$
where -- is the involution of $\nn*\nn$ which
interchanges its two generators. Moreover, these fusion rules were shown
to characterise the $A_u(F)$'s, and this can be interpreted as follows.

\begin{theo}[\cite{cmp}]
The $R^+$-deformations of any $A_u(F)$ are exactly all the $A_u(F)$'s.
\end{theo}

The quotient of $A_u(F)$ by the relations $u=F\bar{u}F^{-1}$ is called
$A_o(F)$. Its universal property shows that it corresponds to the notion of ``algebra of continuous functions on the quantum (or free)
orthogonal group''. See Van Daele-Wang \cite{vdw}. As the operator
$F\bar{F}$ is an intertwiner of the fundamental corepresentation $u$, the
algebra $A_o(F)$ is defined only for matrices $F$ satisfying
$F\bar{F}=$ scalar multiple of the identity (see Wang's paper \cite{wn4} for what happens when $F\bar{F}\notin\c\cdot Id$).

\begin{theo}[\cite{cras}]
The $R^+$-deformations of any $A_o(F)$ are exactly all the $A_o(F)$'s.
\end{theo}

Actually one can easily prove that $A_o\begin{pmatrix}0&1\cr
  -1&0\end{pmatrix}$ is isomorphic to $C(\ss\uu (2))$, so the above result says that the
  $R^+$-deformations of $C(\ss\uu (2))$ are exactly all the $A_o(F)$'s
  (with $n\in\nn$ and $F\in\gg\ll (n,\c)$ satisfying $F\bar{F}\in\c\cdot Id$). For $n=2$ the
  $A_o(F)$'s coincide with Woronowicz's deformations of $C(\ss\uu (2))$ from
  \cite{wo1} (see below). For $n\geq 3$ the Woronowicz algebra
  $A_o(F)$ is a more exotic object: the dimension of its
  irreducible corepresentation $r_k$ corresponding to the
$k$-dimensional irreducible representation of $\ss\uu (2)$ is given by
$$dim(r_k)=\frac{x^k-y^k}{x-y}$$
where $x,y$ are the solutions of $X^2-nX+1=0$. See \cite{cras} and
sections 2 and 5 in \cite{cmp}.

One can ask then whether there
  exist such structure results for the $R^+$-deformations of $\ss\uu
  (N)$ for arbitrary $N$. Besides theorem 2.1 for $G=\ss\uu (N)$,
  which gives examples and theorem 2.3 which solves the problem for $N=2$, we have the following results.

\begin{theo}[\cite{kw}]
Any rigid monoidal semisimple ${\bf C}$-category having the
Grothendieck semiring isomorphic to $\r^+(\ss\uu (N))$ is monoidal
equivalent to a twist $Rep(U_q\sl_N)^\tau$ of the category
$Rep(U_q\sl_N)$ of finite dimensional representations of $U_q\sl_N$, for some $q\in\c^*$ which is not a root of unity, and which is uniquely determined up to $q\leftrightarrow q^{-1}$.
\end{theo}

Here $\tau$ is a $N$-th root of unity (see Kazhdan-Wenzl \cite{kw}). We
mention that for $C(\ss\uu (2))_\mu$ where $\mu\in  [-1,1]-\{ 0\}$ is
as in \cite{wo1} computation of $q$ and $\tau$ gives
$$Corep(C(\ss\uu (2))_\mu )\simeq^{\otimes}
Rep(U_{\sqrt{\mid\mu\mid}}\sl_2)^{sgn(\mu )}$$
For $N\geq 3$ it is not known what values of the
twist can arise from Woronowicz algebras. This should be related to
the question asked by Woronowicz at the end of \cite{wo3}.

\begin{theo}[\cite{sun}]
The $R^+$-deformations of $C(\ss\uu (N))$ are exactly its $R$-matrix
quantizations in the most general sense -- the one of Gurevich \cite{g1}.
\end{theo}

Theorem 2.5 was proved in the following way: theorem 2.4 and the Tannaka-Krein
type duality of Woronowicz \cite{wo3} show that the $R^+$-deformations of
$C(\ss\uu (N))$ are in one-to-one correspondence with the faithful
monoidal functors (satisfying certain positivity conditions) from the
categories $Rep(U_q\sl_N)^\tau$ to the category of finite dimensional
Hilbert spaces. These functors can be shown (via reconstruction) to
be in one-to-one correspondence with the $R$-matrices in \cite{g1} (satisfying certain positivity conditions). In fact the proofs of the above theorems
2.2, 2.3 and of theorem 2.6 below also use reconstruction methods (see section 4 below).

Actually the $R$-matrix quantization of $\ss\uu (N)$ in the spirit
of FRT \cite{frt} requires some work. In \cite{sun} the quantum group
$\ss\uu (N)_R$ was introduced in the same time as an object coming via
duality and (quite sketchy) as a ``compact form'' of Gurevich's $\ss\ll (N,\c )_R$ . See \cite{aa} for
a detailed construction of $\ss\uu (N)_R$. See also \cite{hai} and
\cite{o} for related results, obtained via different methods.

We mention that Gurevich's $R$-matrices (satisfying the same
positivity requirements as those needed in theorem 2.5) were shown
by Wassermann to appear naturally in the theory of full multiplicity
ergodic actions of $\ss\uu (N)$ on von Neumann algebras. A complete
classification of these $R$-matrices for $N=3$ may be found in
\cite{wa2}.

We recall from Wang \cite{wn2} that associated to any finite
dimensional $\c^*$-algebra $B$ is a Woronowicz algebra $A^{aut}(B)$,
which by definition has a universal property making it the
``algebra of countinuous functions on the compact quantum automorphism
group of $B$''. If $n:=dim(B)$ is $1,2,3$ then $B=\c^n$ and
$A^{aut}(B)$ is the algebra of functions on the $n$-th
symmetric group. By \cite{aut} if $n\geq 4$ then $A^{aut}(B)$ is
an $R^+$-deformation of $C(\ss\oo (3))$, and this could be interpreted
in the following way (nothing is lost when stating the less precise
result below, as it is easy to see that $A^{aut}(M_2(\c ))\simeq C(\ss\oo (3))$).

\begin{theo}[\cite{aut}]
For $dim(B)\geq 4$ the $A^{aut}(B)$'s are each other's $R^+$-deformations.
\end{theo}

All these results seem to justify our terminology
``$R^+$-deformation''. However, it is not clear how $R^+$-deformation
could be related to operator algebraic notions of deformation
(see e.g. Rieffel \cite{ri}, Blanchard \cite{bla}, Wang \cite{wn12}). Maybe a
constructive proof of the ``anti-deformation'' conjecture 8.1 below would do part of the job.

\section{Presentation of discrete quantum groups}

In this section we discuss the notion of presentation for finitely
generated (discrete quantum groups represented by) Woronowicz
algebras. We will see in the next section that this kind of
considerations are the first step in the so-called reconstruction
method, which was used for proving theorems 2.2--2.6.

We will use the following notations. If $D$ is a unital $\c$-algebra,
$V$ and $W$ are two finite
dimensional $\c$-linear spaces and $v\in\l (V)\otimes D$ and $w\in\l
(W)\otimes D$ are two elements we define
$$v\otimes w:=v_{13}w_{23}\in \l (V)\otimes\l (W)\otimes D$$
$$Hom(v,w):=\{ T\in \l (V,W)\mid (T\otimes id)v=w(T\otimes id)\}$$
If $D$ is a Hopf $\c$-algebra and $v,w$ are corepresentations, then $\otimes$
and $Hom$ are the usual tensor product, respectively space of
intertwiners. Note also that in general, for any unital $\c$-algebra
$D$, the elements of the form $v\in\l (V)\otimes D$ with $V$ ranging
over the finite dimensional $\c$-linear spaces form a monoidal category
with these $Hom$ and $\otimes$, so our notations are not as abusive as
they seem to be.

We will use freely the terminology from \cite{wo3} concerning concrete
monoidal $\ww^*$-categories, that we will call concrete
monoidal $\c^*$-categories. We recall that the word ``concrete''
comes from the fact that the monoidal $\c^*$-category is given together with
an embedding into (= faithful monoidal $\c^*$-functor to) the
category of finite dimensional Hilbert spaces. We recall also from \cite{wo3} that if $H,K$ are finite dimensional Hilbert
spaces then any invertible antilinear map $j:H\to K$ gives rise to two
linear maps
$$t_j:\c\to H\otimes K,\,\,\, t_j(1)=\sum_i e_i\otimes j(e_i)$$
$$t_{j^{-1}}:\c\to K\otimes H,\,\,\, t_{j^{-1}}(1)=\sum_i f_i\otimes
j^{-1}(f_i)$$
where $\{ e_i\}$, $\{ f_i\}$ are (arbitrary) orthonormal bases in $H$,
respectively $K$. The objects $H$ and $K$ are said to be conjugate in
a concrete monoidal $\c^*$-category $\cat$ (containing
them) if there exists $j:H\to K$ such that
$$t_j\in Hom_\cat(\c ,
H\otimes K)\,\,\,\,\, t_{j^{-1}}\in Hom_\cat (\c ,K\otimes H)$$

The result below is somehow a translation of one of the three
important particular cases of the main result in \cite{wo3} (see section 4 below).

\begin{theo}
Let $n\geq 2$ be an integer, $I$ be a set, $\{ a_k\}_{k\in I}$
and $\{ b_k\}_{k\in I}$ be positive integers, and $X=\{ T_k\mid k\in
I\}$ be a set of
linear maps of the form
$$T_i:(\c^n)^{\otimes a_k}\to (\c^n)^{\otimes b_k}$$
(where $(\c^n)^{\otimes 0}:=\c$). If there exists
$j:\c^n\to\c^n$ such that $t_j,t_{j^{-1}}\in X$ then
$$A_X=\c^*<(u_{ij})_{i,j=1,...,n}\mid u=(u_{ij})\,\,\mbox{is
  unitary},\, T_k\in Hom(u^{\otimes a_k}, u^{\otimes b_k}),\,\forall\,
  k\in I>$$
is a finitely generated Woronowicz algebra. Its category of
  corepresentations is the completion in the sense of
  \cite{wo3} of the concrete monoidal $\c^*$-category $\ca$ defined
  as follows: its objects are the tensor powers of $\c^n$, and its
  arrows are such that $\ca$ is the smallest monoidal category
  containing the arrows $\{ T_k\mid k\in I\}$.

Moreover, any finitely generated Woronowicz algebra arises in
this way.
\end{theo}

The definition of $A_X$ should be understood as follows. Let
$F$ be the free $*$-algebra on $n^2$ variables $(u_{ij})_{i,j=1,...,n}$ and let
$u=(u_{ij})\in\l (\c^n)\otimes F$. By explicitating our notations for
$Hom$ and $\otimes$ with $D =F$ we see that each condition of the form
$T_k\in Hom(u^{\otimes a_k}, u^{\otimes b_k})$, as well as the
condition ``$u$ is unitary'', could be interpreted as being a
collection of relations between the $u_{ij}$'s and their adjoints.  Let
$J\subset F$ be the two-sided $*$-ideal generated by all these
relations. Then the matrix $u=(u_{ij})$ is unitary in $M_n({\bf C}
)\otimes (F/J)$, so its coefficients $u_{ij}$ are of norm  $\leq 1$
for every  ${\bf C}^*$-seminorm on $F/J$ and the enveloping ${\bf
  C}^*$-algebra of $F/J$ is well-defined. We call it $A_X$.

The definition of $\ca$ should be understood as follows: its arrows
are linear combinations of (composable) compositions of tensor
products of maps of the form $T_k$, $T_k^*$ and
$id_m:=$ identity of $(\c^n)^{\otimes m}$. It is clear that $\ca$ is a
concrete monoidal $\c^*$-category.

\begin{proof}
The assumptions $t_j,t_{j^{-1}}\in X$ show that $\c^n$ is a
self-conjugate object in $\ca$. Thus theorem 1.3 in \cite{wo3}
applies, and shows that the universal $\ca$-admissible pair -- which
is $(A_X,u)$ by definition of $A_X$ -- is a finitely generated
Woronowicz algebra (i.e. a compact matrix pseudogroup in the terminology there) whose category of corepresentations is the
completion of $\ca$.

Conversely, let $A$ be an arbitrary finitely generated Woronowicz
algebra. Choose $n\in\nn$ and $v\in M_n(A)$ a unitary corepresentation
whose coefficients generate $A$. By repalcing $v$ with a unitary
representation which is equivalent to the sum $v+\overline{v}$ we may
assume that $v$ is equivalent to $\overline{v}$. Let $X$ be the set of
{\em all} intertwiners between {\em all} tensor powers of $v$. As $v$ is equivalent to $\overline{v}$ we get from
\cite{wo3} the existence of $j:\c^n\to\c^n$ such that
$t_j,t_{j^{-1}}\in X$. Thus $X$ satisfies the condition in the
statement, so we may consider the algebra $A_X$. As both $A$ and $A_X$
can be obtained from $\ca$ by taking the completion and then by
applying duality, we get by uniqueness of these operations (see \cite{wo3}) that they are isomorphic.
\end{proof}

Notice that in the proof of the converse one can take $X$ to be
countable -- just take bases in each space of intertwiners. 

As a first application, we may state the following (non-trivial) definition.

\begin{defi}
A finitely generated Woronowicz algebra is said to be finitely
presented if it is of the form $A_X$ with $X$ a finite set.
\end{defi}

We will see in the next section that the Woronowicz algebras $A_o(F)$,
$A_u(F)$ and $A^{aut}(B)$ are finitely presented, and that the set $X$ has 1
or 2 arrows. In fact, besides these universal quantum groups, many
recently studied ``universal quantum objects'' such as the BMW algebra
\cite{bw} or the Fuss-Catalan algebra \cite{bj1} are ``1- or
2-generated'' in a certain sense (see Bisch-Jones \cite{bj2}).

To our knowledge, the only known examples of non-finitely presented
Woronowicz algebras are those of the form $\c^*(\Gamma )$, with
$\Gamma$ a non-finitely presented discrete group. One can show by using
\cite{ver} and the hereditarity consequence of the main result in
\cite{subf} that Bhattacharyya's planar algebras \cite{bh} have to
come from Woronowicz algebras. The fact that these planar algebras are
not finitely generated (cf. \cite{bh}) should imply that the corresponding
Woronowicz algebras are non-finitely presented.

\section{Reconstruction techniques for computing $R^+$}

In this section we explain how the so-called ``reconstruction method''
for computing $R^+$ of Woronowicz algebras $A$ given with generators and
relations works. This method was used in \cite{wo3} for $C(\ss\uu
(N))_q$, in \cite{cras} for $A_o(F)$, in \cite{cmp} for $A_u(F)$, in
\cite{sun} for $C(\ss\uu (N))_R$ and in \cite{aut} for
$A^{aut}(B)$. Some related techniques were used in \cite{kw} and in
\cite{subf} for proving the main results there. The reconstruction method uses the Tannaka-Krein-Woronowicz duality \cite{wo3} and has four steps.

(I,II) Translate the presentation of $A$ into a ``presentation'' of
its category of corepresentations $Corep(A)$. That is, (I) show that
$Corep(A)$ is the ``smallest'' monoidal category containing certain
arrows (``generators'') and (II) find the relevant formulae (``relations'') satisfied by these arrows.

(III) Brutal combinatorial computation of $Corep(A)$ by using the old method ``use generators and relations for writing
everything as reduced words''.

(IV) In fact what we are interested in is
just the Grothendieck semiring $R^+(A)$ of $Corep(A)$, and in many cases
there exist ad-hoc arguments for ending the computation of
$Corep(A)$ once enough information needed for computing $R^+(A)$ is known.

We will briefly describe how this method is used. There are in fact three cases.

-- Self-adjoint case. (I) We have $A_o(F)=A_X$ where $X$ is the set
consisting of one arrow $E:\c\to\c^n$ given by $1\mapsto\sum
F_{ji}e_i\otimes e_j$. Also $A^{aut}(B)=A_X$ where $X$ is the set
consisting of two arrows: the multiplication $\mu :B\otimes B\to B$ and the
unit $\eta :\c\to B$ (here we use an isomorphism $B\simeq\c^n$ with $n=dim(B)$). (II) The only relevant formula satisfied by $E$ is
$(E^*\otimes id)(id\otimes E)=c\cdot id$, with $c\in\rr$ a
constant. The only relevant formulae satisfied by $\eta$ and $\mu$ are
those coming from the axioms of the algebra structure on $B$. (III)
the relations in (II) and theory from \cite{j1} allow one to show in
both cases that the algebras $End_\ca ((\c^n)^{\otimes k})$ are
Temperley-Lieb-Jones algebras, and that the spaces $Hom_\ca (\c ,(\c^n)^{\otimes
  k})$ have bases indexed by certain non-crossing partitions. (IV) By
(III) we get that the dimensions of the spaces $Hom_\ca (\c ,(\c^n)^{\otimes
  k})$ are given by the Catalan numbers, and an ad-hoc argument (see
below) allows one to end the computation here. See \cite{cras} and \cite{aut}.

-- Non-self-adjoint case. When $A$ is given together with a natural
non-self-adjoint representation, say $u$, it is technically convenient to work
with an analogue of theorem 3.1, where the monoid of objects of
$\ca$ is the free monoid on two copies of $\c^n$ (one for $u$ and
one for $\overline{u}$). This is done in \cite{cmp} for $A_u(F)$: (I) the
set $X$ has two arrows, (II) these arrows satisfy two relations, (III,IV)
what happens here is quite similar to what happens for $A_o(F)$. This
is also done in a quite abstract setting in \cite{subf}. We mention that
the computations in section 4 in there (i.e. step III) were
successfully finished due to the fact that we used in the first stage the
formalism of planar diagrams.

-- Third case. It may happen that $A$ is given together with a natural
non-self-adjoint representation, say $u$, which satisfies $u\subset
u^{\otimes N}$ for some $N$. In this case it is technically convenient to work
with an analogue of theorem 3.1, where the monoid of objects of
$\ca$ is still the one in theorem 3.1, but where complication arises
around complex conjugation of $\c^n$ in $\ca$. This happens in the
``A'' case: for $C(\ss\uu (N))_q$ \cite{wo3} or $U_q\sl_N$ \cite{kw}
or $C(\ss\uu (N))_R$ \cite{sun}. Let us just mention that (I) the set
$X$ consists of one arrow $\c^n\to (\c^n)^{\otimes N}$, (II,III)
the proofs use theory of the Hecke algebra of type $A$ from \cite{f+},
\cite{j12}, \cite{we1}.

We mention that it would be of interest for subfactor theory to say
something about the fusion semiring of Woronowicz algebras coming from
vertex models \cite{ver}, and especially from those coming from spin
models \cite{kac} and from Krishnan-Sunder permutation matrices
\cite{ks}. One can easily see that the ``reconstruction'' method does
not bring anything new, i.e. that the combinatorial problem (III)
coincides with the corresponding combinatorial problem one meets in
subfactor theory (see e.g. \cite{j3}).

We describe now the ``ad-hoc argument'' used as step (IV) in the
reconstruction method for $A_o(F)$, $A_u(F)$ and $A^{aut}(B)$. We recall from \cite{vdn} that if $(P,\phi )$ is a $*$-algebra together with a linear
form, the $*$-distribution of an element $a\in P$ is the functional
$$\mu_a:\c <X,X^*>\,\, \displaystyle{\mathop{\longrightarrow}^{X\mapsto
    a}}\,\, P\,\, \displaystyle{\mathop{\longrightarrow}^\phi}\,\,\c$$
If $(P,\phi )$ is a $\c^*$-algebra together with a faithful state and
    if $a=a^*$ then $\mu_a$ may be viewed (first by restricting it to
    $\c [X]$, then by extending it by continuity) as a probability
    measure on the spectrum of $a$. 

If $(A,u)$ is a finitely generated Woronowicz algebra, one can
consider the $*$-distribution $\mu_{\chi (u)}$ of the character of the
fundamental representation with respect to the Haar functional $\int\in
A^*$. By using the orthogonality formulae in \cite{wo2} we get the following result.

\begin{prop}
The $*$-moments of $\mu_{\chi (u)}$ are given by
$$\mu_{\chi (u)}(M)=\int M(\chi (u),\chi (u)^*)=\int\chi
(M(u,\overline{u}))=dim(Hom(1,M(u,\overline{u})))$$
for any non-commutative 2-variable monomial $M$, where
$M(u,\overline{u})$ denotes the image of $M$ through the unique morphism of
monoids $<X>*<X^*>\to (R^+(A),\otimes )$ given by $X\mapsto u$ and  $X^*\mapsto\overline{u}$.

Thus the pointed semiring $(R^+(A),u)$ uniquely determines $\mu_{\chi (u)}$.\hfill $\Box$
\end{prop}

The interest in $\mu_{\chi (u)}$ comes from the fact that converses
of the last assertion -- which may be used as step (IV) in the
reconstruction method -- hold in certain cases.

\begin{lemm}[cf. \cite{cras,cmp,aut}]
Let $(A,u)$ be a finitely generated Woronowicz algebra.

(i) $\mu_{\chi (u)}$ is semicircular if and only if $R^+(A)\simeq
\r^+(\ss\uu (2))$.

(ii) $\mu_{\chi (u)}$ is circular if and only if $R^+(A)\simeq
R^+(A_u(I_2))$.

(iii) $\mu_{\chi (u)}$ is quarter-circular if and only if $R^+(A)\simeq
\r^+(\ss\oo (3))$.
\end{lemm}

Let us also mention that another point of interest in $\mu_{\chi
  (u)}$ is the following useful lemma, which was the key argument in
  the computation of certain free products of Woronowicz algebras \cite{cmp}, \cite{subf}.

\begin{lemm}[\cite{subf}]
Let $(A,u)$ be a finitely generated Woronowicz algebra and let
$\varphi :A\to B$ be a surjective morphism of Woronowicz algebras. If
$\mu_{\chi (\varphi_*u)}=\mu_{\chi (u)}$ then $\varphi$ is an
isomorphism.
\end{lemm}

The relation of these results with Voiculescu's free probability theory is at a combinatorial level, so it is very unclear. The most conceptual
result in this sense seems to be the one in \cite{subf}, where the operation of going ``from finitely
generated Woronowicz algebras to Popa systems and back'' was
explicitely computed in terms of some free products. This result
also hasn't been understood yet at the spatial (:=
non-combinatorial) level. In fact at the spatial level there is only
the following nice remark, to be related to lemma 4.1 (i): by identifying $\ss\uu (2)$
with the sphere $\ss^3$ we see that the character of the
fundamental representation of $\ss\uu (2)$ is
semicircular.

See Biane \cite{bi} for other applications of free
probability techniques to representation theory.

The explanation for the fact that certain notions related to $(A,u)$ depending
only on $(R^+(A),u)$ turn out to depend only on $\mu_{\chi (u)}$
(cf. e.g. lemmas 4.1 and 4.2 and the comments before theorem 5.2
below) seems to be the fact that knowing $\mu_{\chi (u)}$ imposes
important restrictions on the metric space $Irr(A)$ in proposition 5.1 below.

\section{Applications of $R^+$ to discrete quantum groups}

Proposition 2.1 could be interpreted as saying that any property of a
discrete group $\Gamma$ could be translated in terms of
$R^+(\c^*(\Gamma ))$. One should expect that, more generally, there are
many properties of (discrete quantum groups represented by)
arbitrary Woronowicz algebras $A$ that can be translated in terms of $R^+(A)$.

On the other hand, it is part of the subfactor philosophy that
analytical properties of subfactors should be read from their standard
invariants, and, in certain cases, from their principal graphs or from
their fusion algebras (see e.g. \cite{po3}, \cite{po5}; for fusion algebras of subfactors see e.g. \cite{bi2}). As
subfactors are known to be closely related to Woronowicz algebras (see
e.g. \cite{en}, \cite{i}, \cite{subf}, \cite{ver}, \cite{kac}) this gives real hope for the corresponding properties of
(discrete quantum groups represented by) Woronowicz algebras to
depend only on $R^+$ (or at least only on $(R^+,dim)$, or on $(R^+,
list)$, cf. the comment on principal graphs at the end of section 6 below).

A third reason for believing that this is the case comes from discrete group
philosophy. We recall that if $\Gamma$ is a discrete group of finite
type, any finite set of generators $X\subset\Gamma$ satisfying $1\in
X=X^{-1}$ gives rise to a distance $d_X$ on $\Gamma$. The
quasi-isometry class of the metric space $(\Gamma ,d_X)$ in the sense
of Margulis \cite{ma} does not depend on $X$. A number of properties of
discrete groups were shown to be geometric, i.e. to depend only on
their quasi-isometry class (see e.g. Gromov \cite{gr} and the survey
\cite{gh}). The result below extends the notion of metric and of
quasi-isometry class to the finitely generated (discrete
quantum groups represented by) Woronowicz algebras, in terms of $R^+$.

\begin{prop}
Let $A$ be a finitely generated Woronowicz algebra. Choose a
corepresentation $v\in R^+(A)$ whose coefficients generate $A$, such that $1\subset v=\overline{v}$ (e.g. take
$v=1+u+\overline{u}$ in definition 1.1). Then
$$d_v(a,b)=inf\{ n\in\nn\mid 1\subset a\otimes\overline{b}\otimes
v^{\otimes n}\}$$
makes $Irr(A)$ into a metric space. Moreover, the quasi-isometry class
of $(Irr(A),d_v)$ does not depend on the choice of $v$.
\end{prop}

\begin{proof}
Note first that $d_v(a,b)<\infty$ for any $a,b$ -- just take an
irreducible component $r\subset a\otimes\overline{b}$, and $n\in\nn$
such that $\overline{r}\subset v^{\otimes n}$ (cf. \cite{wo2}). By
Frobenius reciprocity we have
$$d_v(a,b)=inf\{ n\in\nn\mid a\subset v^{\otimes n}\otimes b\} =
inf\{ n\in\nn\mid b\subset v^{\otimes n}\otimes a\}$$
This shows that $d_v(a,b)=d_v(b,a)$ for any $a,b$, and that $d_v(a,b)=0$ if and only if $a=b$. Let $a,b,c\in Irr(A)$. As $b\subset v^{\otimes d_v(a,b)}\otimes a$ and
$b\in v^{\otimes d_v(b,c)}\otimes c$ we get that
$$1\subset b\otimes\overline{b}\subset 
v^{\otimes d_v(a,b)}\otimes a\otimes\overline{c}\otimes  v^{\otimes
  d_v(b,c)}$$
By Frobenius reciprocity we get $1\subset a\otimes\overline{c}\otimes
v^{\otimes (d_v(a,b)+d_v(b,c))}$, and it follows that
$$d_v(a,c)\leq d_v(a,b)+d_v(b,c)$$
Thus $(Irr(A),d_v)$ is a metric space. Let us prove now the last assertion. Let $w\in R^+(A)$ be another
corepresentation whose coefficients generate $A$, such that $1\subset w=\overline{w}$. Let $a,b\in Irr(A)$ be arbitrary elements. We consider the distance
$d_{v+w}(a,b)$ with respect to the corepresentation $v+w\in R^+(A)$ and we have
$$1\subset a\otimes\overline{b}\otimes (v+w)^{\otimes d_{v+w}(a,b)
  }\subset 
a\otimes\overline{b}\otimes (v+v^{\otimes d_v(1,w)})^{\otimes d_{v+w}(a,b)
  }\subset 
a\otimes\overline{b}\otimes
v^{\otimes (1+d_v(1,w))d_{v+w}(a,b) }$$
It follows that $d_v(a,b)\leq Kd_{v+w}(a,b)$ with $K=1+d_v(1,w)$
  independent of $a,b$ and as we clearly have $d_{v+w}(a,b)\leq K^\prime d_v(a,b)$ with
  $K^\prime =1$ we get that $d_v$ and $d_{v+w}$ (hence
$d_v$ and $d_w$ also) are quasi-equivalent.
\end{proof}

If $A=\c^*(\Gamma )$ with $\Gamma$ a discrete group of finite type
then $v$ has to be of the form $\sum_{g\in X}g$, where $X$ is a system of generators of $\Gamma$ satisfying
$1\subset X=X^{-1}$, so the distance $d_v$ on $Irr(A)=\Gamma$ is the usual
distance associated to $X$. When $A=C(G)$ the construction of distances on $\widehat{G}$ is probably well-known, but we were unable to
find a suitable reference; however see McKay \cite{mck} for $\ss\uu
(2)$ and its subgroups. The computation using fusion rules in \cite{cmp} of
$(Irr(A_u(F)),d_{1+u+\overline{u}})$ is left as an exercise to the reader.

We hope that the reader is convinced that many analytical properties of
Woronowicz algebras should depend only on $R^+$ (or at least only on
$(R^+,dim)$ or on $(R^+,list)$; cf. the subfactor point of view, see
the beginning of this section). Unfortunately there seem to be only
two properties for which the translation was done. The first one is
amenability (see section 1 for its definition).

\begin{theo}[G. Skandalis, see \cite{subf}]
Let $A$ be a finitely generated reduced Woronowicz algebra and $u\in M_n(A)$
be a corepresentation whose coefficients generate $A$. Let $\chi
(u)=(Tr\otimes id)u\in A$ be the character of $u$.

(i) The spectrum $X\subset\rr$ of $Re(\chi (u))$ is contained in $[-n,n]$.

(ii) $A$ is amenable if and only if $n\in X$.
\end{theo}

This result is the quantum Kesten theorem. Indeed,
let $\Gamma =<g_1,...,g_n>$ be a finitely generated discrete
group. Then $A=\c^*_{red}(\Gamma )$ is a finitely generated reduced
Woronowicz algebra with fundamental corepresentation $u=diag(\lambda
(g_i))$, where $\lambda :\Gamma\to \c^*_{red}(\Gamma )$ is the left regular
representation. Thus the operator
$Re(\chi (u))$ is exactly the one in Kesten's criterion for the
amenability of $\Gamma$:
$$Re(\chi (u))=\frac{1}{2}\sum_{i=1}^n \lambda (g_i)+\lambda
(g_i^{-1})\in \c^*_{red}(\Gamma )\subset B(l^2(\Gamma ))$$

Let us come back to the general case. The spectrum $X$ of $Re(\chi (u))$ is the support of the spectral
measure $\mu_{Re(\chi (u))}$. As $\mu_{Re(\chi (u))}$ depends only on
the $*$-distribution $\mu_{\chi (u)}$, which in turn depends only on
the pointed semiring $(R^+(A),u)$ (cf. proposition 4.1) we get the following result.

\begin{theo}[\cite{subf}]
Let $A$ be a finitely generated amenable Woronowicz algebra and $B$ be
an $R^+$-deformation of $A$. Then any isomorphism $f:R^+(A)\to R^+(B)$
has to be dimension-increasing, i.e. we have
$$dim(f(u))\geq dim(u)$$
for any $u\in R^+(A)$. Moreover, $B$ is amenable if and only if $f$ is
dimension-preserving.
\end{theo}

We mention that a similar result was independently obtained by Longo
and Roberts at the level of monoidal $\c^*$-categories \cite{lr}.

By combining this with theorems 2.3, 2.1, 2.6 and with the trivial
fact that $C(G)$ is amenable for any compact group $G$ we get the following results.

\begin{coro}[\cite{cmp}, \cite{subf}, \cite{aut}]
(i) $A_o(F)$ is amenable iff $F\in\gg\ll (2,\c )$.

(ii) If $G$ is a simple compact Lie group and $q>0$ then $C(G)_q$ is
amenable.

(iii) If $B$ is a finite dimensional $\c^*$-algebra then $A^{aut}(B)$
is amenable iff $dim(B)\leq 4$.
\end{coro}

A direct proof of (ii) would certainly be quite difficult: for
$G=\ss\uu (N)$ this was done by Nagy in \cite{n} by using direct quite technical arguments.

The second translation concerns Powers' Property of de la Harpe \cite{h}.

\begin{theo}[\cite{cmp}]
Let $A$ be a Woronowicz algebra. We endow the set $\p
(Irr(A))$ of subsets of $Irr(A)$ with the involution $\overline{S}=\{\overline{a}\mid a\in
S\}$ and with the multiplication
$$S\circ T=\{ r\in Irr(A)\mid\exists\, a\in S,\exists\, b\in T\,\,\,
with\,\,\, r\subset a\otimes b\}$$

We say that $A$ has Powers' Property if for any finite subset
$F\subset Irr(A)-\{ 1\}$ there exist elements $r_1,r_2,r_3\in Irr(A)$
and a partition $Irr(A)=D\coprod E$ such that $F\circ D\cap
D=\emptyset$ and $r_s\circ E\cap r_k\circ E=\emptyset$, $\forall\,
s\neq k$.

If $A$ has Powers' Property then $A_{red}$ is simple, with at most one trace.
\end{theo}

Powers' Property for $A$ depends of course only on $R^+(A)$, but this
is not so interesting from the point of view explained in the
beginning of this section, because it is true by definition. However,
theorem 5.3 could be regarded as an illustrating example for the following
general method (with $(P)=$ $A_{red}$ is simple, with at most one trace): 
$$\mbox{compute }R^+(A)\Longrightarrow \mbox{get that
  }A\mbox{ has }(P)$$

Finally, let us remark in connection with discrete groups that any
statement about all compact groups which extends to all compact quantum groups
has to hold for all duals of discrete groups (see section 7 in
\cite{bc} for what may happen in this case). Fortunately, not
all the results on compact groups are of this kind, and it would be interesting
for instance to find the correct quantum extension of the theory of
actions of compact groups on von Neumann algebras from \cite{hls},
\cite{wa1}, \cite{wa12}, \cite{wa2}, \cite{pw}. On the positive (?)
side some results which extend to all compact quantum groups were
given in \cite{bo} and \cite{kac}; on the other hand, some
counterexamples were found in \cite{wn34} and in the first version of
\cite{kac}. We mention that these counterexamples suggest that for certain delicate operator
algebra problems the ``good'' definition for a ``compact quantum
group'' should be (at least) that of a ``co-amenable compact quantum group of Kac type''.

\section{Modular theory, positive parameters and the invariant $(R^+,list)$}

Proposition 2.1 could be interpreted as saying that $R^+$ is a complete
invariant for any $\c^*$-algebra of a discrete group. In this section
we use Woronowicz's work on the modular properties of the Haar
functional and on the square of the antipode \cite{wo2} for
introducing a finer invariant, that we call $(R^+,list)$. We mention that the construction of $(R^+,list)$ is a very simple consequence
of the results in \cite{wo2} and also that this invariant appears (not in a very explicit form) independently in section 5 in Wang's
paper \cite{wn3} and in section 1 in \cite{subf}. 

We first recall from theorem 2.2 that the invariant $R^+$ distinguishes the
$A_u(F)$'s from all the other Woronowicz algebras, but does not
distinguish between the $A_u(F)$'s. The invariant $(R^+,dim)$
distinguishes $A_u(I_2)$ from $A_u(I_3)$. However, it does not distinguish
$A_u(I_2)$ from $A_u\begin{pmatrix}1&0\cr
  0&q\end{pmatrix}$. One invariant which does this job is
$(R^+,qdim)$, where $qdim$ is the quantum dimension of
corepresentations. Conversely, if $q>0$ is
such that $q^2+q^{-2}=3$ then $(R^+,qdim)$ does not distinguish between
$A_u\begin{pmatrix}q&0\cr 0&q^{-1}\end{pmatrix}$ and
$A_u(I_3)$, but these algebras are distinguished by $(R^+,dim)$. The next step is to consider the invariant $(R^+,dim,\,
qdim)$. However, this invariant is not as fine as expected: if $q\in
(0,1)$ and $q_\pm$ is one of the positive solutions of
the equation $X^2+X^{-2}=3\pm q$ then
$$A_u(diag(q_+,q_+^{-1},q_-,q_-^{-1}))$$
have the same $(R^+,dim,\, qdim)$ invariant but are not isomorphic
(see below).

The above results may be obtained by using the
fact that the fusion algebra of $A_u(F)$ is the free $*$-algebra on
its fundamental corepresentation $u$, so both $dim$ and $qdim$ are
uniquely determined by their values on $u$ (cf. \cite{cmp}). See below
for the definition of $qdim$, and for his value $qdim(u)$.

There is one invariant which is finer than both $(R^+,dim)$ and
$(R^+,qdim)$ and that we conjecture it splits the class of Woronowicz
algebras into finte sets. This is the invariant $(R^+,list)$, where
$list(u)$ is the list of eigenvalues of the square root $Q_u$ of the
canonical intertwiner between $u\in R^+$ and its double contragradient. Let us first recall in detail the construction $u\mapsto Q_u$.

\begin{theo}[\cite{wo2}]
Let $A$ be a Woronowicz algebra and let $\a\subset A$ be the dense
$*$-subalgebra of coefficients of finite dimensional corepresentations of $A$. Then there exists a unique family of characters 
$\{ f_z\}_{z\in\c}$ of $\a$ having the following properties (where $*$
denotes the convolutions over the Hopf $\c$-algebra $\a$):

(f1) $f_z*f_{z^\prime}=f_{z+z^\prime}$, $\forall\, z,z^\prime\in\c$
and $f_0=\varepsilon$ (the counit of $\a$).

(f2) the square of the antipode $\kappa$ of $\a$ is given by $\kappa^2(a)=f_{-1}*a*f_1$, $\forall\, a\in\a$.

(f3) $f_z\kappa (a)=f_{-z}(a)$,
$f_z(a^*)=\overline{f_{-\overline{z}}(a)}$, $\forall\, a\in\a$ and
$z\in\c$.

(f4) $\int (ab)=\int (b(f_1*a*f_1))$, $\forall\, a,b\in\a$, where $\int
:A\to\c$ is the Haar functional.
\end{theo}

If $u\in \l (H)\otimes\a$ is a finite dimensional unitary
corepresentation then the restriction of $f_z$ to the
space of coefficients of $u$ can be computed by using the
formula $(id\otimes f_z)v=Q_u^{2z}$, where $Q_u=(id\otimes
f_{\frac{1}{2}})u\in \l (H)$ (this follows from {\it (f1)}). The
following useful characterisation of $Q_u$ follows easily from theorem 6.1.

\begin{lemm}[\cite{subf}, cf. theorem 5.4 in \cite{wo2}]
Let $u\in \l (H)\otimes\a$ be a finite dimensional unitary
corepresentation of $\a$. Then $Q=Q_u=(id\otimes f_{\frac{1}{2}})u$
has the following properties:

(i) $Q>0$ and $Tr(Q^2.)=Tr(Q^{-2}.)$ on $End(u)$.

(ii) $Q^t\,\overline{u}\, (Q^t)^{-1}$ is unitary.

(iii) $Q^2\in Hom(u,(id\otimes\kappa^2)u)$.

\noindent Conversely, if $Q\in\l (H)$ satisfies (i) and one of (ii) and (iii)
then $Q=Q_u$. 
\end{lemm}

This operator $Q_u$ allows one to associate to $u$ some canonical
objects (cf. \cite{subf}). The positive number (cf. (i))
$$qdim(u):=Tr(Q_u^2)=Tr(Q_u^{-2})=(Tr\otimes f_1)u=f_1\chi (u)$$
is called the quantum dimension of $u$. It coincides with the quantum
dimension as defined in \cite{lr}. The linear form (cf. (i))
$$\tau_u:T\mapsto qdim(u)^{-1}Tr(Q_u^2T)=qdim(u)^{-1}Tr(Q_u^{-2}T)$$
is called the canonical trace on $End(u)$. The unitary
corepresentation (cf. (ii))
$$\hat{u}:=(Q_u)^t\,\overline{u}\, (Q_u^t)^{-1}=
(t\otimes f_{\frac{1}{2}}*\kappa (.)*f_{-\frac{1}{2}})u$$
is called the canonical dual of $u$. The list of eigenvalues of $Q_u$
is called the list of positive parameters associated to $u$ and is
denoted $list(u)$ (we call {\em lists} the sets with repetitions,
i.e. a list of $n$ elements of a set $X$ is an element of $X^n/\s_n$).

It is easy to see from theorem 6.1 that
$$ Q_{u+w}=Q_u\oplus Q_w,\, Q_{u\otimes w}=Q_u\otimes Q_w,\, Q_{\hat{v}}=(Q_v^t)^{-1}$$
and it follows that $qdim$, the canonical dual, and $list$ have
some obvious additive, multiplicative and involutive properties (see \cite{subf}). More
precisely, for $list$ we have that  $list(u+w)=list(u)\cup list(w)$
(union of lists), $list(u\otimes w)=list(u)\cdot list(w)$
(multiplication of lists) and $list(\hat{u})=list(u)^{-1}$ (inversion
of lists). It is also clear that $list(u)$ depends only on the
class $u\in R^+(A)$. Summarizing, we have:

\begin{prop}
The map $u\mapsto list (u)$ factors through $R^+(A)$ into a morphism of semirings
$$list:R^+(A)\to (Lists(\rr_+^*),\cdot,\cup )$$
from $R^+(A)$ to the semiring of lists of positive numbers. \hfill $\Box$
\end{prop}

Since $dim(u)=\# list(u)$ the invariant $(R^+,list)$ is finer
than $(R^+,dim)$. Also since
$$qdim(u)=\sum_{q\in list(u)}q^2=\sum_{q\in list(u)}q^{-2}$$
the invariant $(R^+,list)$ is finer than $(R^+,qdim)$. This last remark has the following converse: any list $L$ of positive
numbers satisfying $\sum_{q\in L}q^2=\sum_{q\in L}q^{-2}$ arises as
$list(u)$ for a certain corepresentation $u$ of a certain Woronowicz
algebra $A$ -- just take $u$ to be the fundamental corepresentation of
$A=A_u(Q)$, with $Q$ a diagonal matrix having $L$ as list of
eigenvalues (cf. lemma 6.1). One can get from this and from theorem
2.2 that $(R^+,list)$ is a complete invariant for each $A_u(F)$
(see \cite{wn4}, \cite{subf}). Notice also that $(R^+,list)$ being
finer than $R^+$, it is a complete invariant for each $\c^*(\Gamma )$ with
$\Gamma$ a discrete group (cf. proposition 2.1) and also that it
distinguishes between the $C(G)$'s with $G$ a compact connected Lie group (see
McMullen \cite{mcm} and Handelman \cite{han}). However, $(R^+,list)$ does not distinguish between $C(\ss\uu (2))$ and $C(\ss\uu (2))_{-1}$ (cf. \cite{wo1}).

\begin{conj}[``Finiteness'']
There are only finitely many Woronowicz algebras having a given $(R^+,list)$ invariant.
\end{conj}

We will see in section 8 that this statement is stronger than the
  recently proved result on the finiteness of finite dimensional Kac
  algebras of given dimension (cf. \cite{oc3} and \cite{s}). An even
  finer invariant will be introduced in section 7 below.

We end this section with various considerations on
$(R^+,list)$. First, the definition of $(R^+,list)$ might seem quite
complicated. An equivalent invariant $(R^+,trig)$ may be defined as
follows. By using theorem 6.1 and the additive and multiplicative
properties of $\chi$ we get that the map $trig:R^+(A)\to C(\rr )$ defined by
$$trig(u)(t)= f_{it}(\chi (u))=(Tr\otimes f_{it})u=Tr(Q_u^{2it})=\sum_{q\in
  list(u)}exp(2itq)$$
is a morphism of semirings. The invariants $(R^+,trig)$ and
  $(R^+,list)$ are equivalent.

It follows from theorem 5.2 that the notion of amenability depends only on $(R^+,list)$ (in fact, only on $(R^+,dim)$). Here is another notion encoded in $(R^+,list)$.

\begin{prop}
If $(A,u)$ is a finitely generated Woronowicz algebra then the modular
operator $\Delta_{\int_{A}}$ of the Haar functional $\int_{A}\in A^*$
is diagonal and its point spectrum $Spec(\Delta_{\int_{A}})$ is the
subgroup of $\rr_+^*$ generated by the set $\{ p^2q^2\mid p,q\in list(u)\}$.
\end{prop}

\begin{proof}
Let $l^2(A)$ be the Hilbert space associated to the Haar functional,
and denote by $i:\a\hookrightarrow l^2(A)$ the canonical embedding
(see \cite{wo2}). Let $w\in R^+(A)$ be an arbitrary
corepresentation. By taking a basis in the Hilbert space where $w$
acts consisting of eigenvectors of $Q_w$ we may assume that $w\in
M_n(\a )$ with $Q_w=diag(p_1,..,p_n)$ diagonal. Formula {\it (f4)}
shows that $\Delta_{\int_{A}}$ sends $i(x)\mapsto i(f_1*x*f_1)$ for any $x\in\a$, so its 
restriction to $i(\cat_w)$, where $\cat_w$ is the space of
coefficients of $w$, is given by
$$\Delta_{\int_{A}} :i(w_{ij})\mapsto p_i^2p_j^2i(w_{ij})$$
The Peter-Weyl type decomposition $l^2(A)=\bigoplus_{w\in  Irr(A)}
i(\cat_w)$ (cf. \cite{wo2}) shows that $\Delta_{\int_{A}}$ is diagonal, and that
we have
$$\bigcup_{U\in P}E(U)\subset Spec(\Delta_{\int_{A}})=\bigcup_{w\in Irr(A)}E(w)\subset
\bigcup_{U\in P}E(U)$$
where $E(w)$ denotes the set of eigenvalues of 
$\Delta_{\int_{A}}\mid_{i(\cat_w))}$ for any corepresentation $w$, and where $P$ is the set of corepresentations consisting of all tensor products of $u$'s and $\hat{u}$'s (the last inclusion follows from the fact that 
each $w\in Irr(A)$ is a subcorepresentation of some $U\in P$,
cf. \cite{wo2}). We have shown that $E(w)=\{ p^2q^2\mid p,q\in
list(w)\}$ for any $w$. Together with the additive and multiplicative
properties of $list$ this shows that $E(w\otimes v)=E(w)E(v)$ and $E(\hat{w})=E(w)^{-1}$ (as subsets
of the multiplicative group $\rr_+^*$), so $Spec(\Delta_{\int_{A}})=\bigcup_{U\in
  P}E(U)$ is the group generated by $E(u)$.
\end{proof}

For the Woronowicz algebra $C(\ss\uu (2))_\mu$ from \cite{wo1} we have
$$Spec(\Delta_{\int_{C(\ss\uu (2))_\mu}}) =\mu^{2\zz }$$
(cf. the formulae in the appendix of \cite{wo2}). Together with the
formula $\mu =\tau q^2$ (see the remark after theorem 2.4) this suggests
that for $C(G)_q$ computation should give $Spec(\Delta_{\int_{C(G)_q}})
=q^{c_G\zz}$ with $c_G\in\qq^*_+$ computable. This would be a nice way
of solving the problem ``given $C(G)_q$, what is $q$?''. We
mention that this is a particular case of the more general
problem ``what is a $q$-commutative variety?'' which is currently
being considered (\cite{g2}). Another interesting application is the
following one: if $x_1,...,x_n$ are positive numbers such that $\sum x_k=\sum
x_k^{-1}$, then for $A_u(diag(\sqrt{x_k}))$ we have that
$$Spec(\Delta_{\int_{A_u(diag(\sqrt{x_k}))}})=<\{ x_ix_j\}_{i,j=1...n}>$$
It is not clear how Connes' theory \cite{c} applied to the corresponding von
Neumann algebras could be understood in terms of quantum groups.

We end with a remark on the extra information contained in
$(R^+,list)$. To any $u\in R^+(A)$ one can associate the following tower of
$\c^*$-algebras with traces
$$\c\subset End(u) \subset End(u\otimes\hat{u})\subset
End(u\otimes\hat{u}\otimes u)\subset End(u\otimes\hat{u}\otimes
u\otimes\hat{u})\subset\cdots$$
where the duals and the traces are the canonical ones (see above). By taking the Bratteli diagram and then by deleting the reflections
coming from basic constructions one obtains a weighted graph $\Gamma_u$, called
the principal graph of $u$. It is a standard graph in the sense of subfactor theory
(cf. \cite{subf} and \cite{po4} in the general case; see \cite{kac}
for a constructive proof in the Kac algebra case). It is easy to see
that $\Gamma_u$ depends {\em as a graph} only on the pointed semiring
$(R^+(A),u)$, and {\em as a weighted graph} only on the object $((R^+,list)(A),u)$.

\section{A non-commutative Doplicher-Roberts problem}

We recall that the Woronowicz algebras $C(\ss\uu (2))$ and $C(\ss\uu
(2))_{-1}$ have the same $(R^+,list)$ invariant (cf. \cite{wo1}). They
are not isomorphic: the first one is commutative, and the second one
isn't. A more conceptual argument which shows that they are not
isomorphic is that the ``twist'' in \cite{kw} is given by
$$\tau (Corep(C(\ss\uu (2))_{\pm 1}))=\pm 1$$
(see section 2). That is, $C(\ss\uu (2))$ and $C(\ss\uu (2))_{-1}$ are
distinguished by $Corep$, the monoidal $\c^*$-equivalence class of the
monoidal $\c^*$-category of finite dimensional corepresentations. This invariant $Corep$ is a complete invariant for each $C(G)$ with
$G$ a compact group (cf. Doplicher-Roberts \cite{dr}; see also Deligne
\cite{d1}). It is a also a complete invariant for each $\c^*(\Gamma )$
with $\Gamma$ a discrete group (cf. e.g. proposition 2.1 and the fact
that $Corep$ is finer than $R^+$). However, one can see using $A_u(F)$'s that
$Corep$ is not finer than $(R^+,list)$. Note also that both $Corep$
and $(R^+,list)$ are finer than $(R^+,qdim)$ (see e.g. \cite{lr} for
the definition of $qdim$ in terms of $Corep$).

While the relation between $(R^+,list)$ and $Corep$ is
very unclear, one can consider then the invariant $(Corep,\, list)$, which
is strictly finer than both $(R^+,list)$ and $Corep$. Here $list$
should be viewed as a morphism of semirings from the Grothendieck
semiring of $Corep(.)$ to the semiring of lists of positive
numbers. This is, to our knowledge, the finest known invariant for Woronowicz algebras. Due
to the lack of known examples of Woronowicz algebras, we were unable
to conjecture something on the relation between the invariants
$(Corep,\, list)$ and $Id$. On the other hand, we have the
following much more important and interesting question.

\begin{prob}
What are the possible values of $(Corep,\, list)$?
\end{prob}

This should be regarded as a ``non-commutative Doplicher-Roberts
problem''. Indeed, when restricting attention to algebras of
continuous functions on compact groups, this problem is completely
solved by the Doplicher-Roberts duality (\cite{dr}, see also \cite{dr1}):

-- $Corep$ is a complete invariant, so $list$ is
uniquely determined by it.

-- the values of $Corep$ are all symmetric rigid semisimple
monoidal $\c^*$-categories.

This kind of question is well-known to be quite hopeless for the
moment, due to fact that the examples of semisimple monoidal
$\c^*$-categories which do not come from Woronowicz algebras are many
and varied. Here ``do not come'' means a priori ``do not come by
construction'' but there are tools -- e.g. index or amenability -- for
proving that they really do not come from Woronowicz algebras.

We end this section by recalling what the main examples are. A first kind of
examples come from subfactors of Jones index $< 4$ (see Ocneanu
\cite{oc2}; see also  Kawahigashi \cite{k}). A second kind of examples are
related to Verlinde rules \cite{v} and come from subfactors/quantum groups at
roots of unity (see Wenzl \cite{we1}, \cite{we12}, \cite{we2} and Xu
\cite{xu}) or subfactors/conformal field theory (see Wassermann
\cite{wa3}; see also Jones' Bourbaki expos\'e \cite{jw} and Toledano
\cite{t}). A third kind of examples come from exotic subfactors of
small index $> 4$ (see Asaeda and Haagerup \cite{ah}).

\section{Kac algebras}

Hopf-von Neumann algebras and Hopf $\c^*$-algebras with the square of
the antipode equal to the identity were considered -- long before
quantum groups, Woronowicz algebras, or multiplicative unitaries --
starting with Kac at the beginning of the 60's \cite{ka}. A whole
theory, mainly dedicated to the non-commutative extensions of various
duality results for locally compact abelian groups, was developed since then. See Enock and Schwartz \cite{es} and the references there in.

\begin{defi}
A Woronowicz algebra $A$ is said to be a Woronowicz-Kac algebra if it
satisfies one of the following equivalent conditions (see section 6):

(i) the square of the antipode of $A$ is the identity.

(ii) the Haar functional $\int :A\to\c$ is a trace.

(iii) all $f_z$'s are equal to the counit.

(iv)  the ultraweak closure of $A_{red}$ in its regular representation is a Kac algebra.

(v) $dim(u)=qdim(u)$ for any $u\in R^+(A)$.

(vi) for any $u\in R^+(A)$ the list $list(u)$ consists only of $1$'s
(i.e. ``there are no $q$'s'').
\end{defi}

See Baaj-Skandalis \cite{bs} for the relationship between
Woronowicz-Kac algebras, duals of Woronowicz-Kac algebras, compact Kac
algebras and discrete Kac algebras.

The algebras of continuous functions on compact groups, the
$\c^*$-algebras of discrete groups and the finite dimensional
Woronowicz algebras (= finite dimensional Hopf $\c^*$-algebras) are
Woronowicz-Kac algebras. The universal $\c^*$-algebras $A_u(I_n)$, $A_o(I_n)$ and $A^{aut}(B)$ are also
Woronowicz-Kac algebras. One can construct many Woronowicz-Kac
algebras by using cocycle twistings (see Enock-Vainerman \cite{ev}). A method which produces exotic examples is the use of free
products of discrete quantum groups: for
instance if $G$ and $H$ are compact groups, then the free product
$C(G)*C(H)$ -- which is equal to $\c^*(\widehat{G}*\widehat{H})$ if $G$
and $H$ are commutative -- is a Woronowicz-Kac algebra (see Wang
\cite{wn1}). A big class of examples comes from Jones' vertex models
\cite{j2}: with a suitable terminology, this construction \cite{ver} produces exactly all the algebras representing the ``bi-linear compact
quantum groups of Kac type''. This class of Woronowicz-Kac algebras includes many
exotic objects, such as those associated to vertex models coming from spin models. The Woronowicz algebras
coming from Bhattacharyya's planar algebras \cite{bh} (see section 3) are also of Kac type.

Most of the above examples -- $C(G)$ with $G$ compact Lie group,
$\c^*(\Gamma )$ with $\Gamma$ discrete group of finite type, all the
finite dimensional ones, $A_u(I_n)$, $A_o(I_n)$ and $A^{aut}(B)$,
those coming from vertex models and from planar algebras in \cite{bh}, as well as the (full versions of) free products of such objects --
are finitely generated full Woronowicz-Kac algebras. We thought it
useful to include the list of independent axioms for such objects,
which were of great use in the construction of many of these examples.

\begin{defi}[cf. definitions 1.1 and 8.1]
A finitely generated full Woronowicz-Kac algebra is a unital
$\c^*$-algebra $A$ such that there exists $n\in\nn$ and a unitary matrix $u\in M_n(A)$ satisfying the following conditions:

(i) $A$ is the enveloping $\c^*$-algebra of its $*$-subalgebra
generated by the entries of $u$.

(ii) there exists a $\c^*$-morphism $\delta :A\rightarrow A\otimes_{max} A$ such that $(id\otimes\delta )u=u_{12}u_{13}$.

(iii) there exists a $\c^*$-antimorphism $\kappa :A\to A$ sending
$u_{ij}\leftrightarrow u_{ji}^*$.
\end{defi}

In many operator algebraic situations one may/has to restrict attention from
Woronowicz algebras to Woronowicz-Kac algebras. First, the coactions on
von Neumann algebras of the Woronowicz algebras which are not of Kac
type are quite badly understood (see the comments in the end of
section 5). Second, almost all known constructions of subfactors using
Woronowicz algebras (see \cite{kac} and the references there in) use
in fact Woronowicz-Kac algebras. There are even precise results which assert that for
subfactor constructions one has to restrict attention to
Woronowicz-Kac algebras (theorem 8.5 in \cite{e}, theorem 6.2 in
\cite{subf}, proposition 2.1 in \cite{ver}). There seems to be (a
priori) only one exception: one can associate a Popa system to any
corepresentation of any Woronowicz algebra. However, from Woronowicz
algebras which are not of Kac type one cannot obtain amenable Popa systems (see
\cite{subf}). Also, the subfactors associated to such Popa systems
are not undestood in general in terms of quantum groups (but some
important advances on this subject have been recently made by Ueda \cite{u1}, \cite{u2}).

On the other hand, at the level of known examples (see above) only
$q$-deformations with $q>0$ and the like are not of Kac type. More precisely, to our knowledge, each known example of a Woronowicz algebra is related to a Woronowicz-Kac
algebra. The word ``related'' should be taken in a very vague sense. For instance $A_u(F)$ may be thought (a priori) as being related to
$A_u(I_n)$ just because their presentations look quite the same -- the
identity matrix $I_n$ is just replaced by a matrix ``of parameters''
$F$ (actually the ``minimal'' set of parameters in this example is the
list of eigenvalues of $\sqrt{F^*F}$ modulo multiplication by a scalar, cf. section 6).

\begin{conj}[``Anti-deformation'']
Any Woronowicz algebra is a dimension- preserving $R^+$-deformation of
a Woronowicz-Kac algebra.
\end{conj}

This would be of real interest for certain operator algebraic
problems, in connection with properties which are invariant under
dimension-preserving $R^+$-deformation (see section 5) -- experts know that
life is always much easier in the Kac algebra case, where no
``parameters'' come to complicate the computations. As an example, see the
remarks preceding proposition 8 in \cite{cmp} and the remarks after
theorem A in \cite{subf}. A second point of interest is that the
conjecture would add to the above point of view the fact that ``no
fusion semiring is lost when restricting attention to Woronowicz-Kac
algebras''. A third point of interest is that a constructive proof of this conjecture might establish a link between $R^+$-deformation and operator algebraic notions of deformation. 

Another requirement for a good notion of deformation would be that
Woronowicz-Kac algebras should have some ``rigidity'' properties (weaker of
course than the one of being rigid:=undeformable). This might be the case
for $R^+$-deformation: for Woronowicz -Kac algebras the invariants
$(R^+,list)$ and $(R^+,dim)$ are equivalent, so conjecture 6.1 would imply that any Woronowicz-Kac algebra has
finitely many dimension-preserving $R^+$-deformations in the category
of Woronowicz-Kac algebras.

Notice that this would be stronger than the result on the finiteness of
  finite dimensional Kac algebras of given dimension, which was (a
  version of) one of Kaplansky's longstanding conjectures, and which follows
  from the recent work of Stefan \cite{s} and Ocneanu
  \cite{oc3}. Indeed, it's easy to see that there are finitely many choices
  for the fusion semirings of Kac algebras having a given finite
  dimension.

Finally, let us mention that very little is known about
$R^+$-deformations which are not dimension-preserving. If $A$ is a
Woronowicz algebra, let us call dimension function on $R^+(A)$ any
morphism of monoids $d:R^+(A)\to (\nn ,+,\cdot )$, and standard
dimension function any dimension function which comes from
an $R^+$-deformation. Since the fusion algebra (see section 1) of $C(\ss\uu (2))$ (resp. of $A_u(I_2)$) is
the free algebra (resp. $*$-algebra) on one variable (well-known for
$\ss\uu (2)$; see \cite{cmp} for $A_u(I_2)$) we get that any dimension
function on $R^+(C(\ss\uu (2)))$ (resp. $R^+(A_u(I_2))$) is uniquely
determined by its value on the fundamental corepresentation, so
theorem 2.4 (resp. 2.3) shows that any dimension function is
standard. This kind of statement should be regarded as an accident,
for instance because of the restrictions on standard dimension
functions coming from amenability (theorem 5.2). We mention that for
$C(\ss\uu (N))$ with arbitrary $N$ some other (independent)
restrictions on standard dimension functions should come from
Gurevich's Poincar\'e type duality \cite{g1}.


\begin{thebibliography}{99}

\bibitem{aa} A. Abella and N. Andruskiewitsch, Compact quantum groups
  arising from the FRT construction, preprint.

\bibitem{ah} M. Asaeda and U. Haagerup, Exotic subfactors of finite
  depth with Jones indices $(5+\sqrt{13})/2$ and $(5+\sqrt{17})/2$,
  {\em Commun. Math. Phys.} {\bf 202} (1999), 1-63.

\bibitem{bbs} S. Baaj, E. Blanchard and G. Skandalis, Unitaires
  multiplicatifs en dimension finie et leurs sous-objets, preprint.

\bibitem{bs} S. Baaj and G. Skandalis, Unitaires multiplicatifs et
  dualit\'e pour les produits crois\'es de $C^*$-alg\`ebres, {\em
  Ann. Sci. Ec. Norm. Sup.} {\bf 26} (1993), 425-488.

\bibitem{cras} T. Banica, Th\'eorie des repr\'esentations du
  groupe quantique compact libre $O(n)$, {\em C. R. Acad. Sci. Paris}
  {\bf 322} (1996), 241-244.

\bibitem{cmp}  T. Banica, Le groupe quantique compact libre $U(n)$,
  {\em Comm. Math. Phys.} {\bf 190} (1997), 143-172.

\bibitem{ver} T. Banica, Hopf algebras and subfactors associated to
  vertex models, {\em J. Funct. Anal.} {\bf 159} (1998), 243-266.

\bibitem{subf} T. Banica, Representations of compact
quantum groups and subfactors, {\em J. Reine Angew. Math.} {\bf 509}
(1999), 167-198.

\bibitem{aut} T. Banica, Symmetries of a generic coaction, to appear
  in {\em Math. Ann.}.

\bibitem{sun} T. Banica, A reconstruction result for the $R$-matrix
  quantizations of $SU(N)$, preprint.

\bibitem{kac} T. Banica, Subfactors associated to compact Kac
  algebras, preprint.

\bibitem{bc} P. Baum and A. Connes, $K$-theory for discrete groups,
  {\em London Math. Soc. Lect. Notes} {\bf 135} (1988), 1-20.

\bibitem{bh} B. Bhattacharyya, Krishnan-Sunder subfactors and a new
  countable family of subfactors related to trees, Ph.D. Thesis,
  Berkeley (1998) 

\bibitem{bi} P. Biane, Representations of symmetric groups and free
  probability, {\em Adv. Math.} {\bf 138} (1998), 126-181.

\bibitem{bw} J. Birman and H. Wenzl, Braids, link polynomials and a
  new algebra, {\em Trans. Am. Math. Soc.} {\bf 313} (1989), 249-273.

\bibitem{bi2} D. Bisch, Bimodules, higher relative
commutants and the fusion algebra associated to a subfactor, {\em
  Fields Inst. Commun.} {\bf 13} (1997), 16-63.

\bibitem{bj1} D. Bisch and V. Jones, Algebras associated to intermediate subfactors, {\em Invent. Math.} {\bf 128} (1997), 89-157.

\bibitem{bj2} D. Bisch and V. Jones, Singly generated planar algebras
  of small dimension, to appear in {\em Duke Math. J.}.

\bibitem{bla} E. Blanchard, D\'eformations de $C^*$-alg\`ebres de
  Hopf, {\em Bull. Soc. Math. Fr.} {\bf 124} (1996), 141-215.

\bibitem{bo} F. Boca, Ergodic actions of compact matrix pseudogroups
  on $\c^*$-algebras, {\em Ast\'erisque} {\bf 232} (1995), 93-109.

\bibitem{c} A. Connes, Une classification des facteurs de type III,
  {\em Ann. Sci. Ec. Norm. Sup.} {\bf 6} (1973), 133-252.

\bibitem{d1} P. Deligne, Cat\'egories tannakiennes, in
  ``Grothendieck Festchrift'', Vol. II, Birkhauser (1990), 111-195.

\bibitem{dr1} S. Doplicher and J. Roberts, Endomorphisms of
  $\c^*$-algebras, cross products and duality for compact groups, {\em
  Ann. of Math.} {\bf 130} (1989), 75-119.

\bibitem{dr} S. Doplicher and J. Roberts, A new duality
theory for compact groups, {\em Invent. Math.} {\bf 98} (1989),
157-218.

\bibitem{d} V. Drinfeld, Quantum groups, in Proc. I.C.M. Berkeley
  vol. 1 (1986), 798-820.

\bibitem{er} E. Effros and Z.J. Ruan, Discrete quantum groups. I: The
  Haar measure, {\em Internat. J. Math.} {\bf 5} (1994), 681-723.

\bibitem{e} M. Enock, Inclusions irr\'eductibles de facteurs et
  unitaires multiplicatifs II, {\em J. Funct. Anal.} {\bf 154} (1998), 67-109.

\bibitem{en} M. Enock and R. Nest, Inclusions of factors,
  multiplicative unitaries and Kac algebras, {\em J. Funct. Anal.}
  {\bf 137} (1996), 466-543.

\bibitem{es} M. Enock and J.M. Schwartz, ``Kac algebras and
  duality of locally compact groups'', Springer-Verlag, Berlin (1992)

\bibitem{ev} M. Enock and L. Vainerman, Deformation of a Kac
  algebra by an abelian subgroup, {\em Commun. Math. Phys.} {\bf 178}
  (1996), 571-596.

\bibitem{frt} L. Faddeev, N. Reshetikhin and L. Takhtadzhyan,
  Quantization of Lie groups  and Lie algebras, {\em Leningrad
  Math. J.} {\bf 1} (1990), 193-225.

\bibitem{f+} P. Freyd, D. Yetter, J. Hoste, W. Lickorish, K. Millet
  and A. Ocneanu, A new polynomial invariant of knots and links, {\em
  Bull. Amer. Math. Soc.} {\bf 12} (1985), 239-246.

\bibitem{gh} E. Ghys, Les groups hyperboliques, S\'eminaire Bourbaki
  $n^\circ$722 (1990)

\bibitem{gr} M. Gromov, Infinite groups as geometric objects, in
  Proc. I.C.M. Warszawa vol. 1 (1984), 385-392.

\bibitem{g1} D. Gurevich, Algebraic aspects of the quantum
  Yang-Baxter equation, {\em Leningrad Math. J.} {\bf 2} (1991),
  801-828.

\bibitem{g2} D. Gurevich, Quantum varieties, talk at the conference
  ``Hopf algebras and quantum groups'', Bruxelles (1998)

\bibitem{hai} P. Hai, On matrix quantum groups of type $A_n$,
  preprint.

\bibitem{han} D. Handelman, Representation rings as invariants for
  compact groups and limit ratio theorems for them, {\em
  Internat. J. Math.} {\bf 4} (1993), 59-88.

\bibitem{h} P. de la Harpe, Reduced $\c^*$-algebras of discrete
  groups which are simple with unique trace, {\em Lect. Notes Math.}
  {\bf 1132} (1985), 230-253.

\bibitem{hls} R. H\o egh-Krohn, M. Lanstad and E. St\o rmer, Compact
  ergodic groups of automorphisms, {\em Ann. of Math.} {114} (1981), 75-86.

\bibitem{i} M. Izumi, Goldman's type theorems in index theory, in ``Operator algebras and quantum field theory'', Doplicher, Longo,
  Roberts, Zsido eds., International Press (1997), 249-271.

\bibitem{ji} M. Jimbo, A $q$-difference analog of $U\lie$ and the
  Yang-Baxter equation, {\em Lett. Math. Phys.} {\bf 10} (1985), 63-69.

\bibitem{j1} V. Jones, Index for subfactors, {\em Invent. Math.} {\bf
    72} (1983), 1-25.

\bibitem{j12} V. Jones, Hecke algebra representations of braid groups
  and link polynomials, {\em Ann. of Math.} {\bf 126} (1987), 335-388.

\bibitem{j2}  V. Jones, On knot invariants related to some
  statistical mechanical models, {\em Pacific J. of Math.} {\bf 137}
  (1989), 311-334.

\bibitem{jw} V. Jones, Fusion en alg\`ebres de von Neumann et groupes
de lacets (d'apr\`es A. Wassermann), S\'eminaire Bourbaki $n^\circ$800 (1995)

\bibitem{j3} V. Jones, Planar algebras, I, preprint.

\bibitem{ka} G. I. Kac, Ring groups and the duality principle,
  {\em Proc. Moscow Math. Soc.} {\bf 12} (1963), 259-303.

\bibitem{k} Y. Kawahigashi, On flatness of Ocneanu's connections of
  the Dynkin diagrams and classification of subfactors, {\em
  J. Funct. Anal.} {\bf 127} (1995) 63-107.

\bibitem{kw} D. Kazhdan and H. Wenzl, Reconstructing monoidal
  categories. I.M. Gelfand Seminar,  {\em Adv. in Soviet Math.} {\bf
    16} (1993), 111-136.

\bibitem{ks} U. Krishnan and V.S. Sunder, On biunitary
  permutation matrices and some subfactors of index 9, {\em
  Trans. Amer. Math. Soc.} {\bf 348} (1996), 4691-4736.

\bibitem{kvd} J. Kustermans and A. Van Daele, $\c^*$-algebraic quantum
  groups arising from algebraic quantum groups, {\em
  Internat. J. Math.} {\bf 8} (1997), 1067-1139.

\bibitem{le} S. Levendorskii, Twisted function algebras on a compact
quantum group and their representations, {\em St. Petersburg Math.
J.} {\bf 3} (1992) 405-423.

\bibitem{ls} S. Levendorskii and Y. Soibelman, Algebras of functions
  on compact quantum groups, Schubert cells and quantum tori, {\em
  Commun. Math. Phys.} {\bf 139} (1991), 141-170.

\bibitem{lr} R. Longo and J. Roberts, A theory of dimension, {\em
    K-Theory} {\bf 11} (1997), 103-159.

\bibitem{l} G. Lusztig, Quantum deformations of certain simple modules
  over enveloping algebras, {\em Adv. in Math.} {\bf 70} (1988),
  237-249.

\bibitem{maj} S. Majid, Hopf-von Neumann algebra bicrossproducts, Kac
  algebra bicrossproducts, and the classical Yang-Baxter equations,
  {\em J. Funct. Anal.} {\bf 95} (1991), 291-319.

\bibitem{ma} G.A. Margulis, Groupes discrets d'isom\'etries des
  vari\'et\'es \`a courbure n\'egative, in Proc. I.C.M. Vancouver,
  vol. 2 (1974), 21-34.

\bibitem{mck} J. McKay, Graphs, singularities and finite groups, {\em
    Proc. Symp. Pure Math.} {\bf 37} (1980), 183-186.

\bibitem{mcm} J. McMullen, On the dual object of a compact connected group, {\em Math. Z.} {\bf 185} (1984), 539-552.

\bibitem{n} G. Nagy, On the Haar measure of the quantum
$SU(N)$ group, {\em Commun. Math. Phys.} {\bf 153} (1993), 217-228.

\bibitem{oc2} A. Ocneanu, Quantized groups, string algebras, and
  Galois theory for von Neumann algebras, in {\em Operator
  Algebras and applications 2}, London. Math. Soc. Lect. Notes {\bf
  136} (1988), 119-172.

\bibitem{oc3} A. Ocneanu, There are finitely many paragroups having given
  index and given depth, talk at a conference in Madras (1997).

\bibitem{o} C. Ohn, Quantum $SL(3,\c )$'s with classical
  representation theory, {\em J. Algebra}, to appear.

\bibitem{pwo} P. Podles and S.L. Woronowicz, Quantum deformation of
  Lorentz group, {\em Commun. Math. Phys.} {\bf 130} (1990), 381-431.

\bibitem{po3} S. Popa, Classification of amenable
subfactors of type II, {\em Acta Math.} {\bf 172} (1994), 163-255.

\bibitem{po4} S. Popa, An axiomatization of the lattice of
higher relative commutants of a subfactor, {\em Invent. Math.} {\bf
  120} (1995), 427-445.

\bibitem{po5} S. Popa, Some properties of the symmetric enveloping
  algebra of a subfactor, with applications to amenability and
  property $T$, preprint.

\bibitem{pw} S. Popa and A. Wassermann, Actions of compact Lie
  groups on von Neumann algebras, {\em C. R. Acad. Sci. Paris} {\bf
    315} (1992), 421-426.

\bibitem{ri} M. Rieffel, Continuous fields of $\c^*$-algebras
  coming from group cocycles and actions, {\em Math. Ann.} {\bf 283}
  (1989), 631-643.

\bibitem{r1} M. Rosso, Finite dimensional representations of the quantum
analog of the enveloping algebra of a complex semisimple Lie algebra, {\em Commun. Math. Phys.} {\bf 117} (1988), 581-593.

\bibitem{r2} M. Rosso, Alg\`ebres enveloppantes quantifi\'ees, groupes
quantiques compacts de matrices et calcul differentiel non-commutatif,
{\em Duke Math. J.} {\bf 61} (1990), 11-40.

\bibitem{so} Y. Soibelman, Algebra of functions on a compact quantum group and
its representations, {\em Leningrad Math. J.} {\bf 2} (1991), 161-178.

\bibitem{s} D. Stefan, The set of types of $n$-dimensional
  semisimple and cosemisimple Hopf algebras is finite, {\em
  J. Algebra} {\bf 193} (1997), 571-580.

\bibitem{t} V. Toledano, Fusion of positive energy representations of
  $LSpin_{2n}$, Ph. D. thesis, Cambridge (1997)

\bibitem{u1} Y. Ueda, A minimal action of the compact quantum group
  ${\bf SU}_q (n)$ on a full factor, to appear in {\em
    J. Math. Soc. Japan}.

\bibitem{u2} Y. Ueda, On the fixed point algebra under a minimal free
  product type action of the quantum group ${\bf SU}_q (2)$, preprint.

\bibitem{va} J.M. Vallin, $\c^*$-alg\`ebres de Hopf et
  $\c^*$-alg\`ebres de Kac, {\em Proc. London. Math. Soc.} {\bf 50}
  (1985), 131-174.

\bibitem{vdw} A. Van Daele and S. Wang, Universal quantum
  groups, {\em Internat. J. Math.} {\bf 7} (1996), 255-264.
 
\bibitem{v} E. Verlinde, Fusion rules and modular transformations in
  2D conformal field theory, {\em Nucl. Phys. B} {\bf 300} (1988), 360-376.

\bibitem{vdn} D. Voiculescu, K. Dykema and A. Nica, ``Free
random variables'', CRM Monograph Series $n^{\circ}1$, AMS (1993)

\bibitem{wn1} S. Wang, Free products of compact quantum
groups, {\em Comm. Math. Phys.} {\bf 167} (1995), 671-692.

\bibitem{wn12} S. Wang, Deformations of compact quantum groups via
  Rieffel's quantization, {\em Comm. Math. Phys.} {\bf 178} (1996), 747-764. 

\bibitem{wn2} S. Wang, Quantum symmetry groups of finite spaces,
  {\em Comm. Math. Phys.} {\bf 195} (1998), 195-211.

\bibitem{wn3} S. Wang, Classification of quantum groups $SU_q(n)$,
  {\em J. London Math. Soc.}, to appear.

\bibitem{wn34} S. Wang, Ergodic actions of universal quantum groups on
  operator algebras, to appear in {\em Comm. Math. Phys.}.

\bibitem{wn4} S. Wang, Structure and isomorphic classification of
  compact quantum groups $A_u(Q)$ and $B_u(Q)$, preprint.

\bibitem{wa1}  A. Wassermann, Ergodic actions of compact groups on
  operator algebras I: General theory, {\em Ann. of Math.} {\bf 130}
  (1989), 273-319.

\bibitem{wa12} A. Wassermann, Ergodic actions of compact groups on
  operator algebras II: Classification of full multiplicity ergodic
  actions, {\em Can. J. Math.} {\bf 6} (1988), 1482-1527.

\bibitem{wa2} A. Wassermann, Coactions and Yang-Baxter
equations for ergodic actions and subfactors,  in {\em Operator
  Algebras and applications 2}, London. Math. Soc. Lect. Notes {\bf
  136} (1988), 203-236.

\bibitem{wa3} A. Wassermann, Operator algebras and conformal field
  theory III. Fusion of positive energy representations of $LSU(N)$
  using bounded operators, {\em Invent. math.} {\bf 133} (1998), 467-538.

\bibitem{we1} H. Wenzl, Hecke algebras of type $A_n$ and
subfactors, {\em Invent. Math.} {\bf 92} (1988), 349-383.

\bibitem{we12}  H. Wenzl, Quantum groups and subfactors of
type B,C, and D, {\em Commun. Math. Phys.} {\bf 133} (1990), 383-432.

\bibitem{we2} H. Wenzl, $\c^*$-tensor categories from quantum groups,
  {\em J. Am. Math. Soc.} {\bf 11} (1998), 261-282.

\bibitem{wo1} S.L. Woronowicz, Twisted $SU(2)$ group. An example of a
non-commutative differential calculus. {\em Publ. RIMS Kyoto} {\bf 23}
(1987), 117-181.

\bibitem{wo2} S.L. Woronowicz, Compact matrix
pseudogroups, {\em Commun. Math. Phys.} {\bf 111} (1987), 613-665.

\bibitem{wo3} S.L. Woronowicz, Tannaka-Krein duality
for compact matrix pseudogroups. Twisted $SU(N)$ groups,
{\em Invent. Math.} {\bf 93} (1988), 35-76.

\bibitem{wo4} S.L. Woronowicz, A remark on compact matrix quantum
  groups, {\em Lett. Math. Phys.} {\bf 21} (1991), 35-39.

\bibitem{wo5} S.L. Woronowicz, Compact quantum groups, in
  ``Sym\'etries quantiques'' (Les Houches, 1995), North-Holland,
  Amsterdam (1998), 845-884.

\bibitem{xu} F. Xu, Standard $\lambda$-lattices from quantum groups,
  {\em Invent. Math.} {\bf 134} (1998), 455-487.

\end{thebibliography}
\end{document}